\documentclass{article}
\usepackage{amsmath, amsthm, amssymb, mathtools}
\usepackage{graphicx}
\usepackage{subcaption}
\usepackage[backend=biber,style=alphabetic, doi=false]{biblatex}
\usepackage[hidelinks]{hyperref}
\usepackage{enumitem}
\usepackage{import}
\usepackage[dvipsnames]{xcolor} 
\hypersetup{
    colorlinks, linkcolor=NavyBlue,
    citecolor=OliveGreen, urlcolor=RedOrange
}
\usepackage[nameinlink,capitalise,noabbrev]{cleveref}  
\usepackage[labelfont=bf]{caption}

\newcommand{\subfiglabel}[1]{{\Large \textbf{(#1)}}}

\newsavebox{\subfigstore}
\newcommand{\labeledfig}[2]{%
    \sbox{\subfigstore}{\def\svgwidth{\linewidth}\ignorespaces#2\unskip}%
    \usebox{\subfigstore}%
    \llap{\makebox[\wd\subfigstore][l]{%
        \kern-1.6pt 
        \raisebox{\dimexpr\ht\subfigstore-\height}{\subfiglabel{#1}}%
    }}%
}

\crefname{figure}{Fig.}{Figs.}

\newtheorem{counter}{Counter}[section]

\newtheorem{theorem}[counter]{Theorem}
\newtheorem{lemma}[counter]{Lemma}
\newtheorem{introtheorem}{Theorem}

\newcommand{\RomanNum}[1]{\MakeUppercase{\romannumeral #1}}

\theoremstyle{definition}
\newtheorem{recollection}[counter]{Recollection}
\newtheorem{example}[counter]{Example}
\newtheorem{definition}[counter]{Definition}
\newtheorem{variant}[counter]{Variant}
\newcommand{\WProt}{\operatorname{WProt}}
\newcommand{\Prot}{\operatorname{Prot}}
\newcommand{\Dic}{\operatorname{Dic}}
\newcommand{\D}{\operatorname{D}}
\newcommand{\Zb}{\mathbb{Z}}

\title{Protection Group of Colored Links}

\author{Radin Zahedigolpayegani \and Toni Annala}

\begin{document}
\maketitle

\begin{abstract}
    Colored links provide a mathematical framework for studying the topological stability and classification of three-dimensional topological vortex configurations. Earlier classifications in settings of physical interest relied on system-specific invariants. Here, we introduce a universal invariant of colored links, the \textit{protection group}, interpret earlier classification results as computations of this group, and compute it in new examples.
\end{abstract}

\tableofcontents

\section{Introduction}

Topological vortices are among the most prominent examples of topological defects and have attracted sustained attention in both theoretical and experimental physics for several decades \cite{mermin:1979,anderson-toulouse,salomaa-volovik,ruostekoski-anglin,leanhardt:2003,wu:2025,hall:2026}. Mathematically, a \textit{topological vortex configuration} in three dimensions consists of a \emph{core} $L \subset \mathbb{R}^3$, modeled by a tame link, together with a continuous \textit{order-parameter field}
\[
\Phi \colon \mathbb{R}^3 \setminus L \longrightarrow M,
\]
where the \textit{order-parameter space} $M$ parametrizes the local states of the medium. The interplay between the topology of the vortex core and that of the order-parameter space gives rise to a class of physical structures that are mathematically rich, some of which can be realized experimentally. Correspondingly, the research into this subject remains active to this day \cite{kobayashi:2024, nozaki:2025,abe:2026}.

\begin{figure}[htbp]
    \centering

    \begin{subfigure}{0.45\textwidth}
        \centering
        \labeledfig{a} {%
\begingroup%
  \makeatletter%
  \providecommand\color[2][]{%
    \errmessage{(Inkscape) Color is used for the text in Inkscape, but the package 'color.sty' is not loaded}%
    \renewcommand\color[2][]{}%
  }%
  \providecommand\transparent[1]{%
    \errmessage{(Inkscape) Transparency is used (non-zero) for the text in Inkscape, but the package 'transparent.sty' is not loaded}%
    \renewcommand\transparent[1]{}%
  }%
  \providecommand\rotatebox[2]{#2}%
  \newcommand*\fsize{\dimexpr\f@size pt\relax}%
  \newcommand*\lineheight[1]{\fontsize{\fsize}{#1\fsize}\selectfont}%
  \ifx\svgwidth\undefined%
    \setlength{\unitlength}{250.93984529bp}%
    \ifx\svgscale\undefined%
      \relax%
    \else%
      \setlength{\unitlength}{\unitlength * \real{\svgscale}}%
    \fi%
  \else%
    \setlength{\unitlength}{\svgwidth}%
  \fi%
  \global\let\svgwidth\undefined%
  \global\let\svgscale\undefined%
  \makeatother%
  \begin{picture}(1,0.5535097)%
    \lineheight{1}%
    \setlength\tabcolsep{0pt}%
    \put(0,0){\includegraphics[width=\unitlength,page=1]{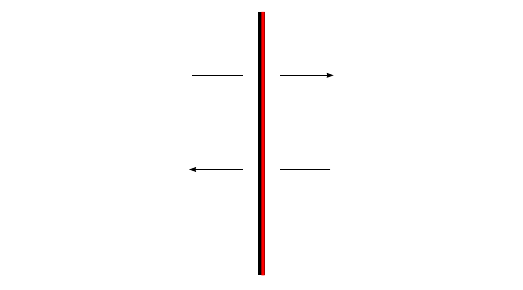}}%
    \put(0.42691066,0.25081178){\color[rgb]{0,0,0}\makebox(0,0)[t]{\lineheight{1.25}\smash{\begin{tabular}[t]{c}$g^{-1}$\end{tabular}}}}%
    \put(0.57299316,0.42955023){\color[rgb]{0,0,0}\makebox(0,0)[t]{\lineheight{1.25}\smash{\begin{tabular}[t]{c}$g$\end{tabular}}}}%
  \end{picture}%
\endgroup%
        }%
        \phantomcaption
        \label{fig:coloring_direction}
    \end{subfigure}
    \hfill
    \begin{subfigure}{0.45\textwidth}
        \centering
        \labeledfig{b} {%
\begingroup%
  \makeatletter%
  \providecommand\color[2][]{%
    \errmessage{(Inkscape) Color is used for the text in Inkscape, but the package 'color.sty' is not loaded}%
    \renewcommand\color[2][]{}%
  }%
  \providecommand\transparent[1]{%
    \errmessage{(Inkscape) Transparency is used (non-zero) for the text in Inkscape, but the package 'transparent.sty' is not loaded}%
    \renewcommand\transparent[1]{}%
  }%
  \providecommand\rotatebox[2]{#2}%
  \newcommand*\fsize{\dimexpr\f@size pt\relax}%
  \newcommand*\lineheight[1]{\fontsize{\fsize}{#1\fsize}\selectfont}%
  \ifx\svgwidth\undefined%
    \setlength{\unitlength}{250.93984529bp}%
    \ifx\svgscale\undefined%
      \relax%
    \else%
      \setlength{\unitlength}{\unitlength * \real{\svgscale}}%
    \fi%
  \else%
    \setlength{\unitlength}{\svgwidth}%
  \fi%
  \global\let\svgwidth\undefined%
  \global\let\svgscale\undefined%
  \makeatother%
  \begin{picture}(1,0.5535097)%
    \lineheight{1}%
    \setlength\tabcolsep{0pt}%
    \put(0,0){\includegraphics[width=\unitlength,page=1]{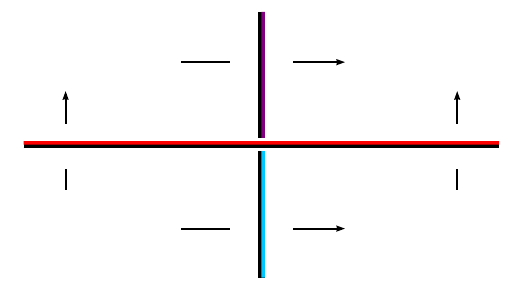}}%
    \put(0.75095393,0.42460534){\color[rgb]{0,0,0}\makebox(0,0)[t]{\lineheight{1.25}\smash{\begin{tabular}[t]{c}$h^{-1}gh$\end{tabular}}}}%
    \put(0.68595432,0.10622431){\color[rgb]{0,0,0}\makebox(0,0)[t]{\lineheight{1.25}\smash{\begin{tabular}[t]{c}$g$\end{tabular}}}}%
    \put(0.15032637,0.37977224){\color[rgb]{0,0,0}\makebox(0,0)[t]{\lineheight{1.25}\smash{\begin{tabular}[t]{c}$h$\end{tabular}}}}%
    \put(0.89943961,0.38237197){\color[rgb]{0,0,0}\makebox(0,0)[t]{\lineheight{1.25}\smash{\begin{tabular}[t]{c}$h$\end{tabular}}}}%
  \end{picture}%
\endgroup%
        }%
        \phantomcaption
        \label{fig:undercrossing_conjugation}
    \end{subfigure}

    \caption{\textbf{Colored link diagrams.}
    (a) The coloring of a strand is a graphical way to record the group element by which an arc in a link diagram is colored. Bicoloring is used to distinguish between a group element and its inverse.
    (b) The coloring has to satisfy Wirtinger-type compatibility condition at crossings: the color of the under-arc changes under conjugation by the color of the arc.}\label{fig:colored_link_diagram}
\end{figure}

\textit{Colored link diagrams} were introduced in \cite{first_paper} as a convenient diagrammatic framework for studying topological vortex configurations. Roughly speaking, a $G$-colored link diagram is a link diagram whose arcs are labeled by elements of a group $G$, subject to compatibility conditions at the crossings that ensure that the coloring records a group homomorphism
\[
\pi_1(\mathbb{R}^3\setminus L) \to G;
\]
see \cref{fig:colored_link_diagram}. In the physical setting, one takes $G=\pi_1(M)$, and the homomorphism is the one induced by the order-parameter field $\Phi$. If $\pi_2(M)=0$, these data determine the homotopy class of $\Phi$, so the colored-link description faithfully represents a topological vortex configuration up to homotopy in this case.

\begin{figure}[htbp]
    \centering
    \captionsetup[subfigure]{skip=0pt} 

    \begin{minipage}[b]{0.48\textwidth}
        \centering
        \offinterlineskip
        \begin{subfigure}{0.9\linewidth}
            \labeledfig{a} {%
                \includegraphics[width=\linewidth]{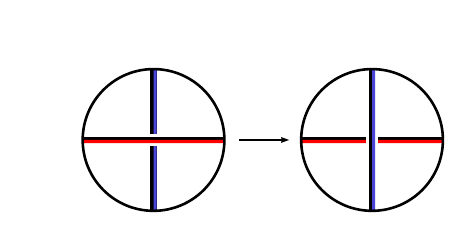}%
            }%
            \phantomcaption\label{fig:strand_crossing_demonstration}
        \end{subfigure}
        \par
        \begin{subfigure}{0.9\linewidth}
            \labeledfig{c} {%
                \includegraphics[width=\linewidth]{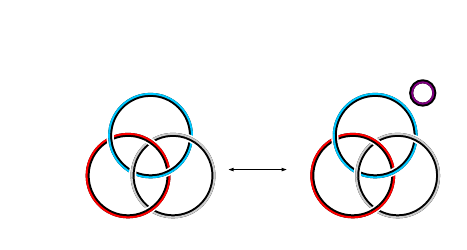}%
            }%
            \phantomcaption\label{fig:unknot_addition}
        \end{subfigure}
    \end{minipage}%
    \hfill
    \begin{minipage}[b]{0.48\textwidth}
        \centering
        \begin{subfigure}{0.9\linewidth}
            \labeledfig{b} {%
                \includegraphics[width=\linewidth]{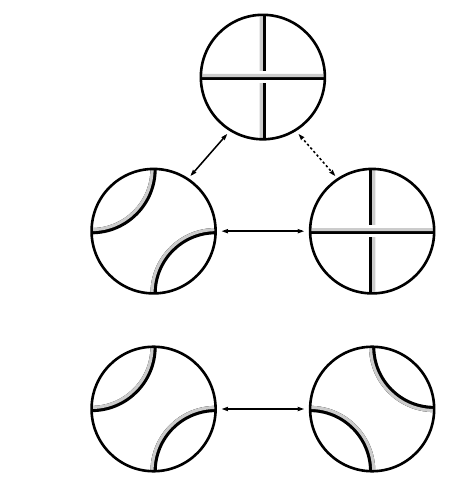}%
            }%
            \phantomcaption\label{fig:local_reconnection_demonstration}
        \end{subfigure}
    \end{minipage}

    \caption{%
        \textbf{Allowed evolution of colored link diagrams.} In addition to the usual knot-theoretic evolution by smooth isotopies and Reidemeister moves, colored link diagrams can evolve in ways that alter the core topology. The two additional moves are (a) local reconnection between strands of the same color provided that the reconnection does not cause a discontinuity of the coloring, and (b) strand crossing when the colors of the crossing strands commute. The elements of the protection group are links up to \emph{stable allowed evolution}, which in addition permits (c) adding and removing simple unlinked loops.%
    }\label{fig:stably_allowed_evolution}
\end{figure}

In the same work, the notions of \textit{topological protection} and \textit{topological equivalence} were introduced to describe vortex evolution that is permitted by the topology of the physical system. A vortex configuration is topologically protected if it cannot evolve into a disjoint union of simple unlinked vortex loops through such allowed evolution. Two configurations are topologically equivalent if one can be transformed into the other by allowed vortex evolution (\cref{fig:stably_allowed_evolution}).  

These notions have since been studied by constructing invariants of colored links that are conserved under allowed vortex evolution \cite{first_paper,second_paper,third_paper}. The purpose of this paper is to construct the universal invariant of this kind, which we call the \emph{protection group} of a symmetric quandle $Q$ that specifies the coloring scheme. Abstractly, $\Prot(Q)$ is the commutative monoid whose elements are equivalence classes of $Q$-colored link diagrams, up to allowed evolution of $Q$-colored link diagrams and addition and removal of unlinked simple loops, and whose monoid operation is given by disjoint union of diagrams. By our first result, it is the universal abelian-group-valued invariant of $Q$-colored links.

\begin{introtheorem}
The monoid $\Prot(Q)$ is an abelian group. Furthermore, any abelian-group-valued invariant of $Q$-colored link diagrams
\[
i \colon Q\text{-}\mathrm{Diagram} \longrightarrow G
\]
that is conserved under the allowed evolution of $Q$-colored link diagrams and that is additive under disjoint unions factors uniquely through a homomorphism $\overline i \colon \Prot(Q) \to G$.
\end{introtheorem}

In addition to defining the protection group, we also compute it in several examples.

\begin{introtheorem}
We have the following isomorphisms:
\begin{enumerate}
    \item odd dihedral groups: $\Prot(\D_{2k + 1}) \cong 0$;
    \item odd dicyclic groups: $\Prot(\Dic_{2k + 1}) \cong 0$;
    \item dihedral groups of degree 2 mod 4: $\Prot(\D_{4k+2}) \cong (\Zb/2)^{2k}$.
\end{enumerate}
\end{introtheorem}

In the language developed here, the classification results of the earlier papers on the subject can be interpreted as computations (sometimes partial) of the corresponding protection groups. Indeed,
\begin{enumerate}
    \item \cite{first_paper} shows that $\Prot(Q_8) \cong \Zb/4$;
    \item \cite{second_paper} shows that $\Prot(\mathrm{tricoloring}) \cong \Zb/3$, where $\mathrm{tricoloring}$ is the symmetric quandle given by the subquandle of transpositions in $\Sigma_3$; and \item \cite{third_paper} shows that $\Prot(T^*_\mathrm{black})$ is generated by the classes of the $T^*$-colored trefoil and the figure-eight knot but does not specify the relations. Here, $T^*_\mathrm{black}$ is the symmetric quandle given by the subset of the binary tetrahedral group $T^*$ that consists of elements of the form $\tfrac{1 \pm i \pm j \pm k}{2}$.
\end{enumerate}
We expect many more computations of the protection group to be carried out in the future.

\subsection*{Declaration of AI use.} Artificial intelligence tools were used to polish the exposition and to turn hand-drawn figures into editable vector-graphics files. For the latter purpose, GPT-6 pro was used.

\section{Quandles and Colored Links}
Here, we recall the notions of a quandle and a symmetric quandle (which is,
in rough terms, a quandle equipped with an inverse operation) and define the 
notion of a $Q$-coloring. 

A large portion of this section is taken from \cite{carter_quandles}; 
specifically, we use Carter's definition of quandles, symmetric quandles, 
quandle coloring of links and coloring of links with symmetric quandles.
\subsection{Quandles and Symmetric Quandles}

A quandle is an algebraic structure that abstracts the notion of a group equipped with its conjugation operation.
\begin{definition}
    A \emph{quandle} $(Q,\triangleleft)$ is a set $Q$ equipped with a binary operation $\triangleleft \colon Q \times Q \to Q$ satisfying the following properties:
    \begin{enumerate}
        \item $a \triangleleft a = a$ for all $a \in Q$;
        \item for all $a, b \in Q$, there exists a unique $c \in Q$ such that $c \triangleleft a = b$;
        \item for all $a, b, c \in Q$:
            \begin{equation*}
                (a \triangleleft b) \triangleleft c = (a \triangleleft c) \triangleleft (b \triangleleft c).
            \end{equation*}
    \end{enumerate}
\end{definition}

\begin{example}
A group $G$ gives rise to a quandle by setting $g \triangleleft h = h^{-1}gh.$
\end{example}

One can define the operation $\triangleleft^{-1}$ by setting $a \triangleleft^{-1} b = c$ if and only if $a = c \triangleleft b$. This is well-defined because $\bullet \triangleleft b$ is bijective.

We now define a symmetric quandle. 
A symmetric quandle, as per \cite{carter_quandles}, is roughly a quandle with a notion of an inverse element.

\begin{definition}
A \emph{symmetric quandle} $(Q,\triangleleft, \rho)$ is a quandle equipped with an involution of sets $\rho \colon Q \to Q$ such that $\rho(a) \triangleleft b = \rho(a \triangleleft b)$ and $a \triangleleft \rho(b) = a \triangleleft^{-1} b$. Such a $\rho$ is referred to as a \emph{good involution}.
\end{definition}

We will usually omit the conjugation (and the good involution, if the quandle is symmetric) operation from the notation when easily understood from context, and merely denote the quandle by $Q$. 

\subsection{\texorpdfstring{$Q$}{Q}-colored link diagrams}

Recall that a tame link can be represented by a link diagram with finitely many arcs and crossings. Next, we define what a $Q$-coloring of such a diagram is, where $Q$ is a symmetric quandle. This definition is taken from \cite{carter_quandles}.

\begin{figure}[htbp]
    \centering

    \begin{subfigure}{0.45\textwidth}
        \centering 
        \labeledfig{a} {%
\begingroup%
  \makeatletter%
  \providecommand\color[2][]{%
    \errmessage{(Inkscape) Color is used for the text in Inkscape, but the package 'color.sty' is not loaded}%
    \renewcommand\color[2][]{}%
  }%
  \providecommand\transparent[1]{%
    \errmessage{(Inkscape) Transparency is used (non-zero) for the text in Inkscape, but the package 'transparent.sty' is not loaded}%
    \renewcommand\transparent[1]{}%
  }%
  \providecommand\rotatebox[2]{#2}%
  \newcommand*\fsize{\dimexpr\f@size pt\relax}%
  \newcommand*\lineheight[1]{\fontsize{\fsize}{#1\fsize}\selectfont}%
  \ifx\svgwidth\undefined%
    \setlength{\unitlength}{170.07874016bp}%
    \ifx\svgscale\undefined%
      \relax%
    \else%
      \setlength{\unitlength}{\unitlength * \real{\svgscale}}%
    \fi%
  \else%
    \setlength{\unitlength}{\svgwidth}%
  \fi%
  \global\let\svgwidth\undefined%
  \global\let\svgscale\undefined%
  \makeatother%
  \begin{picture}(1,0.81666667)%
    \lineheight{1}%
    \setlength\tabcolsep{0pt}%
    \put(0,0){\includegraphics[width=\unitlength,page=1]{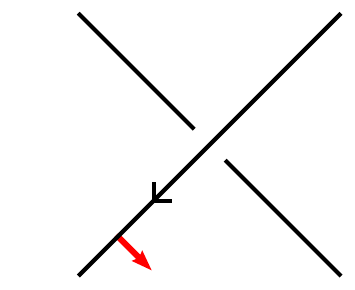}}%
    \put(0.84074594,0.23474636){\color[rgb]{0,0,0}\makebox(0,0)[t]{\lineheight{1.25}\smash{\begin{tabular}[t]{c}$a \triangleleft b$\end{tabular}}}}%
    \put(0.824228,0.53701445){\color[rgb]{0,0,0}\makebox(0,0)[t]{\lineheight{1.25}\smash{\begin{tabular}[t]{c}$b$\end{tabular}}}}%
    \put(0.45512913,0.61156195){\color[rgb]{0,0,0}\makebox(0,0)[t]{\lineheight{1.25}\smash{\begin{tabular}[t]{c}$a$\end{tabular}}}}%
  \end{picture}%
\endgroup%
        }%
        \phantomcaption
        \label{fig:formal_directed_undercrossing_conjugation}
    \end{subfigure}
    \hfill
    \begin{subfigure}{0.45\textwidth}
        \centering
        \labeledfig{b} {%
            \includegraphics[width=\linewidth]{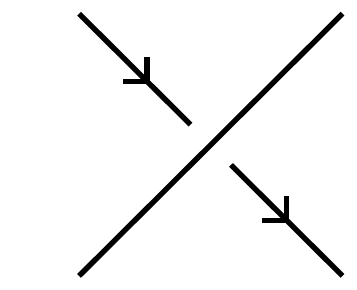}%
        }%
        \phantomcaption
        \label{fig:consistent_orientation}
    \end{subfigure}
    
    \caption{(a) Conjugation of $a$ upon passing under $b$. Note that the normal vector is drawn with respect to the canonical orientation of the plane, i.e., rotated by $90^\circ$ counterclockwise. (b) An orientation is consistent if it does not change upon passing under another strand.}
\end{figure}

\begin{definition}
    Let $K$ be a consistent orientation given to a diagram (\cref{fig:consistent_orientation}). 
    Then a normal vector to $K$ can be drawn at each point of the arcs in the canonical orientation (\cref{fig:formal_directed_undercrossing_conjugation}).

    Let $\mathcal{R}$ be the set of arcs. A directed $Q$-coloring $\mathcal{C} \colon \mathcal{R} \to Q$ is a map such that for each crossing, 
    the color of the undercrossing arc is conjugated upon passing under the 
    other strand in the direction of the normal vector (\cref{fig:formal_directed_undercrossing_conjugation}).
\end{definition}


The above definition requires specifying an orientation, which is more data than we want to specify. This can be fixed by allowing the orientation to flip, as long as the corresponding label is inverted using the involution $\rho$ of a symmetric quandle. 

We formalize the change of direction for a component.  Suppose $Q$ is a symmetric quandle for the following arguments.

\begin{definition}
    Let $K$ be a consistent orientation given to a diagram and $\mathcal{C}$ a directed coloring thereof.
    Let $r_{1}, r_{2}, \dots, r_{k}$ be the different arcs into which a component is broken. We define the operation of \emph{flipping} by changing the direction of this component entirely and changing $\mathcal{C}(r_i)$ to $\rho(\mathcal{C}(r_i))$ for each arc.
\end{definition}

It is clear that \emph{flipping} a component yields a valid directed coloring, as $a = c \triangleleft \rho(b)$ in \cref{fig:formal_directed_undercrossing_conjugation}.

We now define $Q$-colorings of a link for a symmetric quandle $Q$.

\begin{definition}
    Let $(K_1, \mathcal{C}_1)$ and $(K_2, \mathcal{C}_2)$ be two directed $Q$-colorings.
    We define them to be \emph{equivalent} if one can be obtained from the other by flipping a sequence of components.
    It is evident that this relation is an equivalence relation.
    Each resulting equivalence class is called a \emph{$Q$-coloring} of the diagram.
\end{definition}

\begin{figure}[htbp]
    \centering
    \captionsetup[subfigure]{skip=0pt}

    \begin{subfigure}{\dimexpr \textwidth * 100 / 141 \relax}
        \centering
        \labeledfig{a} {%
\begingroup%
  \makeatletter%
  \providecommand\color[2][]{%
    \errmessage{(Inkscape) Color is used for the text in Inkscape, but the package 'color.sty' is not loaded}%
    \renewcommand\color[2][]{}%
  }%
  \providecommand\transparent[1]{%
    \errmessage{(Inkscape) Transparency is used (non-zero) for the text in Inkscape, but the package 'transparent.sty' is not loaded}%
    \renewcommand\transparent[1]{}%
  }%
  \providecommand\rotatebox[2]{#2}%
  \newcommand*\fsize{\dimexpr\f@size pt\relax}%
  \newcommand*\lineheight[1]{\fontsize{\fsize}{#1\fsize}\selectfont}%
  \ifx\svgwidth\undefined%
    \setlength{\unitlength}{283.46456693bp}%
    \ifx\svgscale\undefined%
      \relax%
    \else%
      \setlength{\unitlength}{\unitlength * \real{\svgscale}}%
    \fi%
  \else%
    \setlength{\unitlength}{\svgwidth}%
  \fi%
  \global\let\svgwidth\undefined%
  \global\let\svgscale\undefined%
  \makeatother%
  \begin{picture}(1,0.31)%
    \lineheight{1}%
    \setlength\tabcolsep{0pt}%
    \put(0,0){\includegraphics[width=\unitlength,page=1]{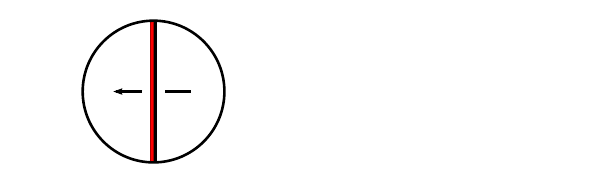}}%
    \put(0.22944694,0.16913716){\color[rgb]{0,0,0}\makebox(0,0)[t]{\lineheight{1.25}\smash{\begin{tabular}[t]{c}$i$\end{tabular}}}}%
    \put(0,0){\includegraphics[width=\unitlength,page=2]{coloring_guide_for_primary_q8.pdf}}%
    \put(0.51456695,0.17406365){\color[rgb]{0,0,0}\makebox(0,0)[t]{\lineheight{1.25}\smash{\begin{tabular}[t]{c}$j$\end{tabular}}}}%
    \put(0,0){\includegraphics[width=\unitlength,page=3]{coloring_guide_for_primary_q8.pdf}}%
    \put(0.80982696,0.16913716){\color[rgb]{0,0,0}\makebox(0,0)[t]{\lineheight{1.25}\smash{\begin{tabular}[t]{c}$k$\end{tabular}}}}%
  \end{picture}%
\endgroup%
        }%
        \phantomcaption
        \label{fig:coloring_guide_for_primary_q8}
    \end{subfigure}%
    \begin{subfigure}{\dimexpr \textwidth * 41 / 141 \relax}
        \centering
        \labeledfig{b} {%
\begingroup%
  \makeatletter%
  \providecommand\color[2][]{%
    \errmessage{(Inkscape) Color is used for the text in Inkscape, but the package 'color.sty' is not loaded}%
    \renewcommand\color[2][]{}%
  }%
  \providecommand\transparent[1]{%
    \errmessage{(Inkscape) Transparency is used (non-zero) for the text in Inkscape, but the package 'transparent.sty' is not loaded}%
    \renewcommand\transparent[1]{}%
  }%
  \providecommand\rotatebox[2]{#2}%
  \newcommand*\fsize{\dimexpr\f@size pt\relax}%
  \newcommand*\lineheight[1]{\fontsize{\fsize}{#1\fsize}\selectfont}%
  \ifx\svgwidth\undefined%
    \setlength{\unitlength}{116.22047244bp}%
    \ifx\svgscale\undefined%
      \relax%
    \else%
      \setlength{\unitlength}{\unitlength * \real{\svgscale}}%
    \fi%
  \else%
    \setlength{\unitlength}{\svgwidth}%
  \fi%
  \global\let\svgwidth\undefined%
  \global\let\svgscale\undefined%
  \makeatother%
  \begin{picture}(1,0.75609756)%
    \lineheight{1}%
    \setlength\tabcolsep{0pt}%
    \put(0,0){\includegraphics[width=\unitlength,page=1]{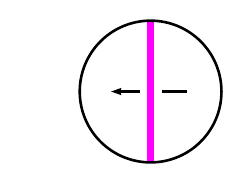}}%
    \put(0.52074864,0.4125332){\color[rgb]{0,0,0}\makebox(0,0)[t]{\lineheight{1.25}\smash{\begin{tabular}[t]{c}$-1$\end{tabular}}}}%
  \end{picture}%
\endgroup%
        }%
        \phantomcaption
        \label{fig:bicoloring_is_single}
    \end{subfigure}

    \begin{subfigure}{\textwidth}
        \centering
        \labeledfig{c} {%
            \includegraphics[width=\linewidth]{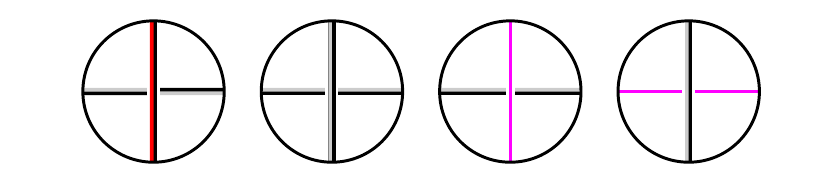}
        }%
        \phantomcaption
        \label{fig:crossing_guide_for_q8}
    \end{subfigure}

    \caption{
    \textbf{$Q_8$-coloring.}
    (a) Coloring for $i$, $j$, and $k$. Note that the color corresponding to undercossing to the other diraction is given by the inverse element. The strand is $i$ in the black-to-red direction and $-i$ in the red-to-black direction.
    (b) $\rho(-1) = -1$, making a single color sufficient.
    (c) Strands cross according to the rule above. Note that the coloring for $i$ flips upon passing under $j$; this is because $i \triangleleft j = -i$, and so must be colored in the other direction.
    }\label{fig:Q8_coloring}
\end{figure}

\begin{example}[$Q_8$-colorings]
The coloring scheme used in the case of the quaternion group $Q_8 = \{\pm 1, \pm i, \pm j, \pm k\}$ is recalled in \cref{fig:Q8_coloring}. 
\end{example}

\begin{figure}[htbp]
    \centering
    \includegraphics[width = \textwidth]{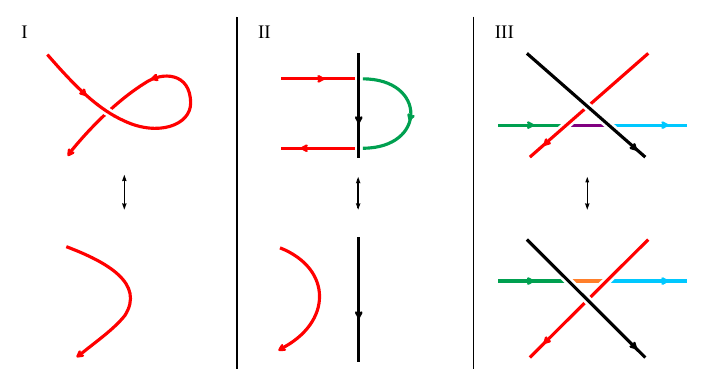}
    \caption{\textbf{$Q$-colored Reidemeister moves.}}
    \label{fig:directed_coloring_under_reidemeister_moves}
\end{figure}

\subsection{Allowed Evolution of \texorpdfstring{$Q$}{Q}-Colored Link Diagrams}
Here, we define the allowed evolution of $Q$-colored links. The evolutions we intend to formalize are of the following four types:
\begin{itemize}
    \item planar isotopies;
    \item Reidemeister moves;
    \item strand crossings;
    \item local reconnections.
\end{itemize}
$Q$-colored links behave trivially under planar isotopies, and Reidemeister moves are explained in \cref{fig:directed_coloring_under_reidemeister_moves}. We limit our focus to the latter two. 

Strand crossings can occur when two crossing strands are colored by commuting elements. However, we do not immediately have a suitable notion of commutativity in quandles. The notion can be generalized quite simply.

\begin{definition}
    Let $Q$ be a quandle, and let $a, b \in Q$. 
    We say that $a$ and $b$ \emph{commute} if and only if $a \triangleleft b = a$ and $b \triangleleft a = b$.
\end{definition}

This generalizes the standard notion of commutativity. 
Note that commutativity is definitionally symmetric.

\begin{lemma}
    Let $a, b \in Q$ commute. Then $\rho(a) \triangleleft b = \rho(a)$, 
    $a \triangleleft \rho(b) = a$, and $\rho(a) \triangleleft \rho(b) = \rho(a)$.
\end{lemma}
\begin{proof}
    \begin{align*}
        \rho(a) &= \rho(a \triangleleft b) \\
                &= \rho(a) \triangleleft b.
    \end{align*}

    Moreover, $a \triangleleft b = a$ implies that $a = a \triangleleft^{-1} b = a \triangleleft \rho(b)$.
    $\rho(a) \triangleleft \rho(b) = \rho(a)$ now follows from the two previous identities.
\end{proof}

This proves that if one flips a component, commutativity of the strands does not change.
Therefore, commutativity is well-defined for $Q$-colorings.

Strand crossing can now be easily defined.
\begin{definition}
    Let $L$ be a $Q$-colored link such that there exists a disk in which two commuting strands cross orthogonally (as in \cref{fig:strand_crossing_demonstration} and \cref{fig:local_reconnection_demonstration}). Let $(K, \mathcal{C})$ be a directed coloring representing the coloring of the link. Let the overcrossing strand have the color $c_1$ and the undercrossing strands have the color $c_2$ with respect to $\mathcal{C}$. Then the $Q$-colored link $L'$ is obtained from $L$ by changing the crossing in the disk, maintaining the orientation $K$, and coloring the previously undercrossing strand by $c_2$ and the two new undercrossing arcs by $c_1$.
\end{definition}

\begin{definition}
    Let $L$ be a $Q$-colored link diagram. Local reconnection is by reconnecting two crossing strands of the same color as in \cref{fig:local_reconnection_demonstration} and extending the former coloring to the new one afterwards.
\end{definition}
    
\subsection{Topological and Stable Topological Equivalence}
The introduced operations allow vortices to evolve over time. 
Two $Q$-colored links are said to be equivalent 
if one can be transformed into the other via the allowed operations. 
We already have a good notion of equivalence classes when the link changes by planar isotopies or Reidemeister moves (i.e., when the link changes under homotopy), 
whereas this class might change under the rest of the evolutions. We address this by introducing the following definition.

\begin{definition}
    Two $Q$-colored, planar tame diagrams are said to be \emph{immediately equivalent} if they are related by a planar isotopy, a Reidemeister move, a strand crossing or a local reconnection.
\end{definition}

Therefore, a link is immediately equivalent to another link if it is the result of exactly one valid operation. It is clear that this relation is reflexive and symmetric. However, since we allow multiple transformation steps, we must allow a link $L_{2}$ to be equivalent to $L_{1}$ if it is the result of multiple transformations applied to $L_{1}$. Therefore, we consider topological equivalence as the one immediate equivalence induces.

\begin{definition}
    Let $L$ and $L'$ be two $Q$-colored, planar tame links. We say $L$ and $L'$ are \emph{topologically equivalent}
    if there is a sequence of $Q$-colored links $L_0, L_1, \dots, L_k$ such that $L_0 = L$, $L_k = L'$, and $L_{i - 1}$ is immediately equivalent to $L_{i}$ for each $1 \leq i \leq k$.
\end{definition}

Stable topological equivalence permits a wider range of transformations by incorporating the addition and removal of arbitrary isolated unknots. 

\begin{definition}
    Let $L$ and $L'$ be two $Q$-colored, planar tame links. We say $L$ and $L'$ are \emph{stably equivalent} and write $L \cong L'$ if there is a sequence of topological equivalence classes $C_0, C_1, \dots, C_k$, with $L \in C_0$ and $L' \in C_k$, such that for each $1 \le i \le k$ there is an element of $C_i$ obtained from an element of $C_{i - 1}$ through the addition or removal of an isolated unknot.
\end{definition}

\section{Protection Group}

\noindent Two $Q$-colored link diagrams can be \emph{added} by putting them alongside each other. This provides an algebraic structure on the set of stable topological equivalence classes of colored links.

\begin{definition}
Let $Q$ be a symmetric quandle. We define $\Prot(Q)$ as the commutative monoid of stable topological equivalence classes of $Q$-colored link diagrams, equipped with the addition operation induced by disjoint union. 
\end{definition}

\begin{figure}[htbp]
    \centering
    \includegraphics[width = 0.5\textwidth]{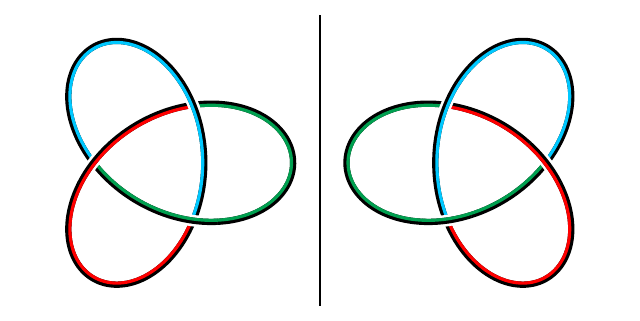}
    \caption{Mirroring illustration}
    \label{fig:mirroring_illustration}
\end{figure}

Our next goal is to show that $\Prot(Q)$ is an abelian group. We can \emph{mirror} a colored link diagram by reflecting it entirely across an arbitrary line in the plane. An example demonstrating this operation is provided in \cref{fig:mirroring_illustration}. 

\begin{theorem}
\label{thm:mirror_is_inverse}
The mirror of a $Q$-colored link is its additive inverse in $\Prot(Q)$.
\end{theorem}

\begin{figure}[htbp]
    \centering

    \begin{subfigure}{0.45\textwidth}
        \centering
        \labeledfig{a} {%
\begingroup%
  \makeatletter%
  \providecommand\color[2][]{%
    \errmessage{(Inkscape) Color is used for the text in Inkscape, but the package 'color.sty' is not loaded}%
    \renewcommand\color[2][]{}%
  }%
  \providecommand\transparent[1]{%
    \errmessage{(Inkscape) Transparency is used (non-zero) for the text in Inkscape, but the package 'transparent.sty' is not loaded}%
    \renewcommand\transparent[1]{}%
  }%
  \providecommand\rotatebox[2]{#2}%
  \newcommand*\fsize{\dimexpr\f@size pt\relax}%
  \newcommand*\lineheight[1]{\fontsize{\fsize}{#1\fsize}\selectfont}%
  \ifx\svgwidth\undefined%
    \setlength{\unitlength}{311.81102362bp}%
    \ifx\svgscale\undefined%
      \relax%
    \else%
      \setlength{\unitlength}{\unitlength * \real{\svgscale}}%
    \fi%
  \else%
    \setlength{\unitlength}{\svgwidth}%
  \fi%
  \global\let\svgwidth\undefined%
  \global\let\svgscale\undefined%
  \makeatother%
  \begin{picture}(1,0.95454545)%
    \lineheight{1}%
    \setlength\tabcolsep{0pt}%
    \put(0,0){\includegraphics[width=\unitlength,page=1]{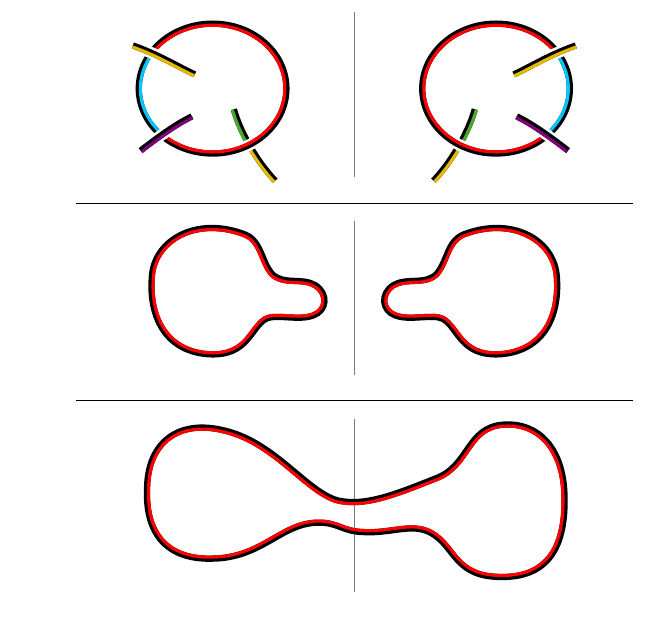}}%
    \put(0.56028814,0.01190652){\color[rgb]{0,0,0}\makebox(0,0)[t]{\lineheight{1.25}\smash{\begin{tabular}[t]{c}$\ell$\end{tabular}}}}%
    \put(0.56028814,0.65008834){\color[rgb]{0,0,0}\makebox(0,0)[t]{\lineheight{1.25}\smash{\begin{tabular}[t]{c}$\ell$\end{tabular}}}}%
    \put(0.56028814,0.34554288){\color[rgb]{0,0,0}\makebox(0,0)[t]{\lineheight{1.25}\smash{\begin{tabular}[t]{c}$\ell$\end{tabular}}}}%
  \end{picture}%
\endgroup%
        }%
        \phantomcaption
        \label{fig:connecting_components_to_l}
    \end{subfigure}
    \hfill
    \begin{subfigure}{0.45\textwidth}
        \centering
        \labeledfig{b} {%
\begingroup%
  \makeatletter%
  \providecommand\color[2][]{%
    \errmessage{(Inkscape) Color is used for the text in Inkscape, but the package 'color.sty' is not loaded}%
    \renewcommand\color[2][]{}%
  }%
  \providecommand\transparent[1]{%
    \errmessage{(Inkscape) Transparency is used (non-zero) for the text in Inkscape, but the package 'transparent.sty' is not loaded}%
    \renewcommand\transparent[1]{}%
  }%
  \providecommand\rotatebox[2]{#2}%
  \newcommand*\fsize{\dimexpr\f@size pt\relax}%
  \newcommand*\lineheight[1]{\fontsize{\fsize}{#1\fsize}\selectfont}%
  \ifx\svgwidth\undefined%
    \setlength{\unitlength}{311.81102362bp}%
    \ifx\svgscale\undefined%
      \relax%
    \else%
      \setlength{\unitlength}{\unitlength * \real{\svgscale}}%
    \fi%
  \else%
    \setlength{\unitlength}{\svgwidth}%
  \fi%
  \global\let\svgwidth\undefined%
  \global\let\svgscale\undefined%
  \makeatother%
  \begin{picture}(1,0.95454545)%
    \lineheight{1}%
    \setlength\tabcolsep{0pt}%
    \put(0,0){\includegraphics[width=\unitlength,page=1]{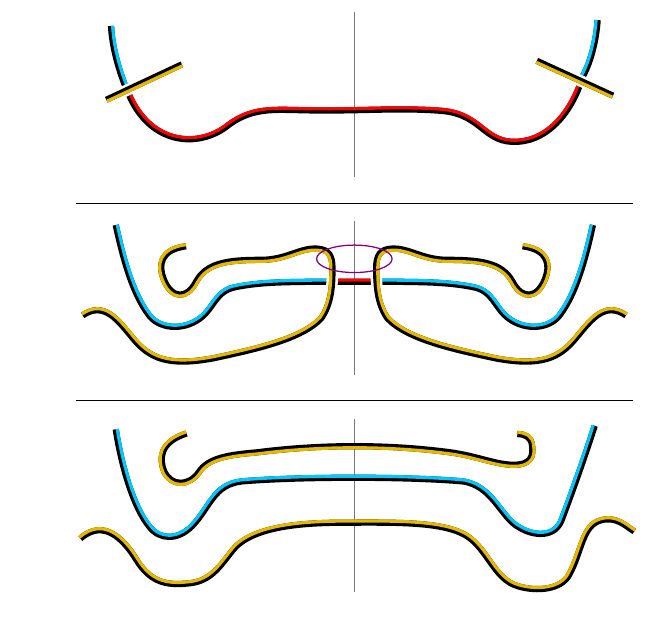}}%
    \put(0.56028814,0.65008834){\color[rgb]{0,0,0}\makebox(0,0)[t]{\lineheight{1.25}\smash{\begin{tabular}[t]{c}$\ell$\end{tabular}}}}%
    \put(0.56028814,0.34554288){\color[rgb]{0,0,0}\makebox(0,0)[t]{\lineheight{1.25}\smash{\begin{tabular}[t]{c}$\ell$\end{tabular}}}}%
    \put(0.56028814,0.01190652){\color[rgb]{0,0,0}\makebox(0,0)[t]{\lineheight{1.25}\smash{\begin{tabular}[t]{c}$\ell$\end{tabular}}}}%
  \end{picture}%
\endgroup%
        }%
        \phantomcaption
        \label{fig:resolving_crossing_near_l}
    \end{subfigure}

    \caption{(a) One can drag the components and attach them to their mirror image by Reidemeister moves and local reconnections.
    (b) Deformation of a component toward to line followed by applications of local reconnections and strand crossings. Note that the format of the new components has not changed; each component on the left still seems \emph{attached} to the line $\ell$.}
\end{figure}

\begin{proof}
    Choose a line $\ell$ such that the diagram lies entirely on one side of it. Mirror the diagram with respect to this line. We proceed by performing Reidemeister moves on the original link and simultaneously repeating them on the mirror to preserve symmetry, and applying local reconnections to strands near the line $\ell$.

    First pick a point from each component of the original link and drag it near the line $\ell$. Apply local reconnections after every such action to attach the components together (\cref{fig:connecting_components_to_l}). In this manner, one can consider each component of the original link as a strand having both of its endpoints connected to the line $\ell$. We will henceforth only speak of the right side of line $\ell$ and implicitly assume the components are connected to the line, as the left side is uniquely determined by symmetry and this symmetry shall be preserved. In our figures, however, we demonstrate both sides for further clarity.

    We proceed by eliminating crossings one by one. If no crossings remain, the problem is solved and the result is the trivial class. Otherwise, consider a component which participates in a crossing. Consider one of its endpoints and traverse the component to the other endpoint until the first crossing is encountered. Without loss of generality, one can assume it is an undercrossing. Then one can drag this crossing along the component to bring it near the line $\ell$ (\cref{fig:resolving_crossing_near_l}), where a local reconnection and a Reidemeister move eliminate the crossing. The symmetry is preserved, and one can thus continue this procedure until no crossings remain, and the resulting link is trivial.
\end{proof}

\begin{figure}
    \centering
    \includegraphics[width=\linewidth]{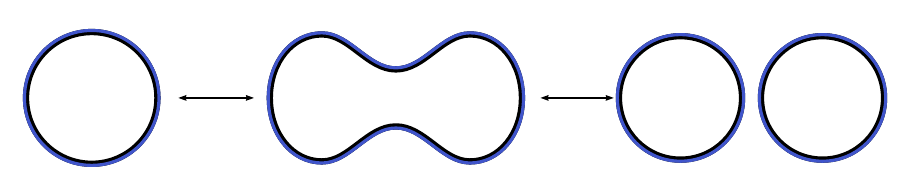}
    \caption{An unknot can be transformed into two through local reconnection.}
    \label{fig:unknot_is_zero}
\end{figure}

\begin{theorem}
    Let $\varphi: Q\text{-}\operatorname{Diagram} \to G$ be a group-valued invariant on $Q$-colored diagrams which respects topological equivalence of links and is also additive on disjoint unions. Then this invariant factors uniquely through $\iota: Q\text{-}\operatorname{Diagram} \to \Prot(Q)$ as a group homomorphism.
\end{theorem}
\begin{proof}
    If $L_1$ and $L_2$ are two diagrams for which $\iota(L_1) = \iota(L_2)$, then by definition they belong to the same stable topological equivalence class, meaning they can be transformed into each other through a series of evolutions and additions and removals of isolated unknots. $\varphi$ remains invariant under evolution; we shall prove it remains invariant under addition and removal of unknots as well.

    Let $U$ be an unknot. This unknot can evolve into two unknots by local reconnection (\cref{fig:unknot_is_zero}). Therefore, $\varphi(U) + \varphi(U) = \varphi(U)$, meaning $\varphi(U) = 0 \in G$. Therefore, $\varphi$ is zero on all unknots, proving their addition and removal preserves the value of $\varphi$.

    Hence $\varphi(L_1) = \varphi(L_2)$, which implies $\varphi$ factors through $\iota$ uniquely as a set map $\overline{\varphi}$. This set map is furthermore a group homomorphism. This is because if $g_1, g_2 \in \Prot(Q)$ and $g_1 = \iota(L_1)$ and $g_2 = \iota(L_2)$ for some $Q$-colored diagrams $L_1, L_2$, then: 
    \begin{align*}
        \overline{\varphi}(g_1 + g_2) &= \overline{\varphi}(\iota(L_1 \sqcup L_2)) \\
        &= \varphi(L_1 \sqcup L_2) \\
        &= \varphi(L_1) + \varphi(L_2) \\
        &= \overline{\varphi}(\iota(L_1)) + \overline{\varphi}(\iota(L_2)) \\
        &= \overline{\varphi}(g_1) + \overline{\varphi}(g_2),
    \end{align*}
    which proves the map is indeed a group homomorphism. 
\end{proof}

\begin{variant}[Weak protection group]
A variant of the protection group is obtained by allowing local reconnections but not strand crossings. Since the proof of \cref{thm:mirror_is_inverse} uses only local reconnections, it shows that the \emph{weak protection monoid} $\WProt(Q)$, consisting of equivalence classes of $Q$-colored links under this finer equivalence relation, is in fact a group. We call it the \emph{weak protection group}. It is often useful to consider $\WProt(Q)$ in place of $\Prot(Q)$. For example, the cobordism invariants defined in \cite{second_paper} are invariants of the weak protection group. In many cases of interest, including tricoloring, one has $\WProt(Q)=\Prot(Q).$
\end{variant}

\section{Computations}
In this section, we compute $\operatorname{Prot}(Q)$ when $Q$ is $\operatorname{Dic}_{2k + 1}$, $\D_{2k + 1}$, or $\D_{4k + 2}$. All of our computations are based on techniques we shall call \emph{surgery} (first introduced in \cite{first_paper}) and \emph{switching}, which switches a (suitable) crossing at the expense of adding an unknot. This technique is loosely introduced in the second subsection, and we compute the protection group for the named quandles in the subsections that follow.

\subsection{Recollection of the Dihedral and Dicyclic Groups}
We work with the following presentation for the dihedral group:

\begin{definition}
    The \emph{dihedral group} of order $2n$, denoted $\D_n$, is given by the presentation

    \[ \D_n = \langle r, s \mid r^n = s^2 = 1,\, srs = r^{-1} \rangle. \]
\end{definition}

$\D_n$ can be realized as the group of symmetries of a regular $n$-gon in the plane.
Such a polygon has $2n$ symmetries: the identity symmetry, $n - 1$ non-trivial rotations, and $n$ reflections.
Rotations are represented by elements of the form $r^k$ in the group (with $r$ being the rotation by $\frac{2\pi}{n}$) and reflections are elements of the form $sr^k$ (with $s$ being a fixed reflection), all for $0 \le k < n$. Clearly, $r^0$ is the identity.

\begin{figure}[htbp]
    \centering
    \includegraphics[width = \textwidth]{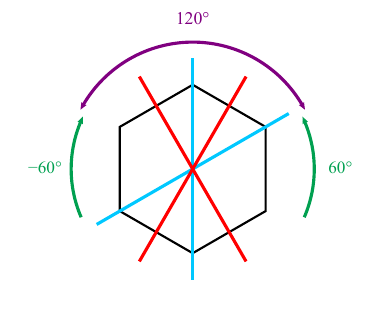}
    \caption{Some elements of $\D_6$. Elements of the same color belong to the same conjugacy class.}
    \label{fig:dihedral_geometric_representation}
\end{figure}

The \emph{dicyclic group} $\operatorname{Dic}_n$ of order $4n$ is defined by the following presentation.

\begin{definition}
    $\operatorname{Dic}_n = \langle a, x \mid a^n = x^2, x^4 = 1, xax^{-1} = a^{-1} \rangle$
\end{definition}

The dicyclic group can then be realized as a subgroup of the quaternions $\mathbb{H}$ with $a = e^{\frac{\pi i}{n}}$, a primitive $2n$-th root of unity, and $x = j$. In particular, $Q_8 \cong \operatorname{Dic}_2$.

In our notation, $\operatorname{Dic}_n$ has elements of the form $a^k, xa^k, -a^k, -xa^k$ for $0 \le k < n$. Here, $-1$ denotes the central involution $a^n = x^2$ in $\operatorname{Dic}_n$.

\subsection{Switching and Surgery}
As switching and surgery are merely proof techniques and not rigorous theorems, 
we settle for intuitive definitions. 

\subsubsection{Switching}
We say two arcs are of the same \emph{type} if one's color 
(in either direction) is the conjugate to the other (in either direction). 
Clearly, different arcs of one component are all of the same type.

\begin{figure}[htbp]
    \centering
    \def\svgwidth{\textwidth}
    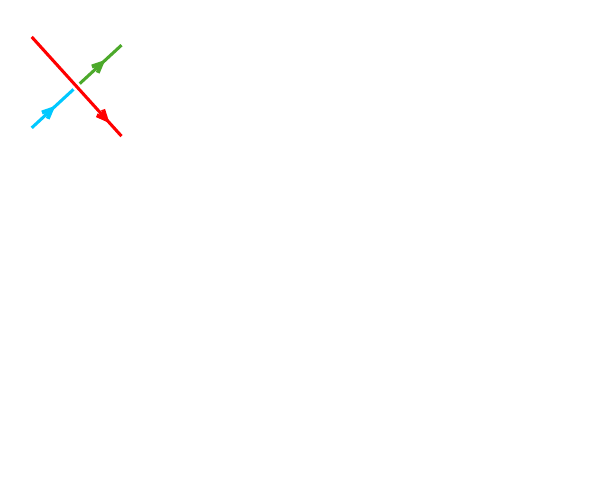
    \caption{The crossing can be switched to either form. 
    Note that $c$ is automatically converted to $a$ as 
    $a = a \triangleleft (b \triangleleft x) = a \triangleleft a$.}
    \label{fig:switching_guide}
\end{figure}

Suppose a crossing is formed by an arc of (directed) color $a$ passing over arcs of colors $b$ and $c$ such that $a = b \triangleleft x$ for some $x \in Q$. 
One can add an unknot colored by $x$ and slide it through the crossing to conjugate $b, c$ into $a$. Then, this crossing can be \emph{switched} to different states as all strands are of the same color. This process is demonstrated in \cref{fig:switching_guide}. Note that switching can happen in two different ways, called type~\RomanNum{1} and type~\RomanNum{2} switching.

Switching can be used to connect different components of a particular conjugacy class of colors and to change crossings in order to unknot a component. This technique will be extensively used in the following subsections.

\subsubsection{Surgery}
\begin{figure}[htbp]
    \centering
    \includegraphics[width = \textwidth]{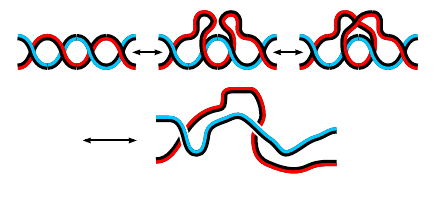}
    \caption{In a $Q_8$-colored link, the $i$-colored strand winds around the $k$-colored strand twice. As a result, its color returns to its initial, and the windings can be removed by surgery.}
    \label{fig:surgery_guide}

\end{figure}

If a strand winds around an arc twice such that its (directed) color changes to its inverse after each undercrossing, \emph{surgery} can be applied to eliminate these windings (\cref{fig:surgery_guide}). 
This technique was introduced in \cite{first_paper}, and we shall use it extensively in the subsequent computations.

\subsection{Computation of \texorpdfstring{$\Prot(\D_{2k+1})$}{Prot(D(2k+1))}}
\begin{recollection}
    The conjugacy classes of $\D_{2k + 1}$ are given by:

    \begin{itemize}
        \item $\{ 1 \}$
        \item $\{ r^t, r^{2k + 1 - t} \}$ for $1 \le t \le k$
        \item $\{ sr^t \mid 0 \le t \le 2k \}$
    \end{itemize}
\end{recollection}

We now have the necessary tool set to compute $\Prot(\D_{2k + 1})$.

\begin{theorem} 
    \label{thm:prot_for_d_2k+1}
    $\Prot(\D_{2k + 1}) \cong 0$, the trivial group.
\end{theorem}

For each $\D_n$, we call the components whose arcs are colored by elements of the form $sr^k$ \emph{primary} (this property is evidently independent of the chosen directed coloring). Note that for each component, either all or none of the arcs are colored by a primary element.

Let $L$ be an arbitrary $\D_{2k + 1}$-colored link. 
We prove the following lemmas.

\begin{figure}[htbp]
    \centering
    \def\svgwidth{\textwidth}
\begingroup%
  \makeatletter%
  \providecommand\color[2][]{%
    \errmessage{(Inkscape) Color is used for the text in Inkscape, but the package 'color.sty' is not loaded}%
    \renewcommand\color[2][]{}%
  }%
  \providecommand\transparent[1]{%
    \errmessage{(Inkscape) Transparency is used (non-zero) for the text in Inkscape, but the package 'transparent.sty' is not loaded}%
    \renewcommand\transparent[1]{}%
  }%
  \providecommand\rotatebox[2]{#2}%
  \newcommand*\fsize{\dimexpr\f@size pt\relax}%
  \newcommand*\lineheight[1]{\fontsize{\fsize}{#1\fsize}\selectfont}%
  \ifx\svgwidth\undefined%
    \setlength{\unitlength}{328.81889764bp}%
    \ifx\svgscale\undefined%
      \relax%
    \else%
      \setlength{\unitlength}{\unitlength * \real{\svgscale}}%
    \fi%
  \else%
    \setlength{\unitlength}{\svgwidth}%
  \fi%
  \global\let\svgwidth\undefined%
  \global\let\svgscale\undefined%
  \makeatother%
  \begin{picture}(1,0.31896552)%
    \lineheight{1}%
    \setlength\tabcolsep{0pt}%
    \put(0,0){\includegraphics[width=\unitlength,page=1]{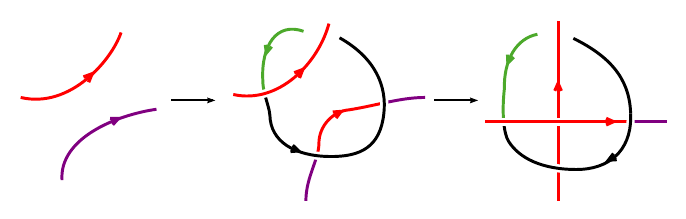}}%
    \put(0.12660742,0.1641335){\color[rgb]{0,0,0}\makebox(0,0)[t]{\lineheight{1.25}\smash{\begin{tabular}[t]{c}$sr$\end{tabular}}}}%
    \put(0.16495833,0.10555302){\color[rgb]{0,0,0}\makebox(0,0)[t]{\lineheight{1.25}\smash{\begin{tabular}[t]{c}$sr^3$\end{tabular}}}}%
    \put(0.57618476,0.08913677){\color[rgb]{0,0,0}\makebox(0,0)[t]{\lineheight{1.25}\smash{\begin{tabular}[t]{c}$r^6$\end{tabular}}}}%
  \end{picture}%
\endgroup%

    \caption{An example of joining two components in $\D_7$}
    \label{fig:joining_primary_components}
\end{figure}

\begin{lemma}
    The link $L$ can evolve into a link $L'$ containing at most one primary component.
    \label{lmm:joining_primary_components}
\end{lemma}
\begin{proof}
    The claim is obvious if no such component exists. Suppose otherwise.

    Two components can be unified by first placing part of one over the other and applying type~\RomanNum{2} switching at the resulted crossing.
    Since the components are colored by elements from the same conjugacy class, one can slide an unknot with the color $r^k$ such that the two arcs attain the same color, and then join the two by a local reconnection (\cref{fig:joining_primary_components}). By repeating this operation, all primary components are joined, and only one primary component remains.
\end{proof}

Apply \cref{lmm:joining_primary_components} so that at most one primary component remains.
If there are no such components, the case is trivial: all other components are colored by rotational elements that commute, and thus can be dragged out and removed by strand crossings and removal of isolated unknots.
Assume, therefore, that exactly one such component exists.

\begin{lemma}
    The primary component can be unknotted.
    \label{lmm:trivializing_the_primary}
\end{lemma}
\begin{proof}
    It is known that each component can be made an unknot by changing its crossings finitely many times.
    This process can be replicated by type~\RomanNum{1} switching the crossings using auxiliary unknots colored by rotational elements. 
    Therefore, one can make the component an unknot with finitely many switches. 
\end{proof}

We will henceforth deform the primary component into an oval for pictorial simplicity. Each component of color $r^k$ winds a number of times around the primary component, changing color to $r^{-k}$ each time in any arbitrary directed coloring. Suppose that $k \ne 0$, as otherwise the component commutes with every element and can be dragged out and removed.
If it winds around the primary component $l \ge 2$ times, surgery can reduce it to $l - 2$ (other components do not play a role as they commute with $r^k$). 
Therefore, for every such component, $l$ can be assumed to be $1$ or $0$. 
It cannot, however, be $1$, as otherwise the coloring of the component would not match with itself if $k \ne 0$\footnote{
    This idea is taken from \cite{first_paper}.
}
Therefore, each $l$ is $0$ and the primary component is entirely unlinked, making it removable. 
After that, the link is solely a collection of components of the form $r^k$, which can be unlinked and eliminated as before. This completes the proof.

\subsection{Computation of \texorpdfstring{$\Prot(\Dic_{2k+1})$}{Prot(Dic(2k+1))}}
\begin{recollection}
    The conjugacy classes of $\operatorname{Dic}_{2k + 1}$ are:

    \begin{itemize}
        \item $\{ 1 \}$
        \item $\{ -1 \}$
        \item $\{ a^t, -a^{2k + 1 - t} \}$ for $1 \le t \le 2k$
        \item $\{ xa^{2t} \mid 0 \le t \le 2k \}$
        \item $\{ xa^{2t + 1} \mid 0 \le t \le 2k \}$
    \end{itemize}

    where $-1$ corresponds to $a^{2k + 1}$. Note that it commutes with every other element.
\end{recollection}

\begin{theorem}
    $\Prot(\operatorname{Dic}_{2k + 1}) \cong 0$.
\end{theorem}

The proof is analogous to that of $\D_{2k + 1}$; we therefore merely give the general outline. We designate as \emph{primary} those components colored by elements of the form $xa^k$. One can similarly turn all primary components into one primary unknot \footnote{
    Note that all arcs with (directed) colors of the class $\{xa^{2k + 1}\}$ 
    can be considered of the class $\{xa^{2k}\}$, if considered in the other direction.
}
(the problem is trivial if no such component exists) and then,
all components not of color $1$ or $-1$ must wind around the primary component evenly many times, switching from $a^k$ to $a^{-k}$ with each pass. This number can be reduced to zero with surgery, transforming the initial link to a trivial state. By these methods, the protection group can be computed in the same manner as for $\D_{2k + 1}$.

\subsection{Computation of \texorpdfstring{$\Prot(\D_{4k+2})$}{Prot(D(4k+2))}}
\begin{recollection}
    The conjugacy classes of $\D_{4k + 2}$ are given by:

    \begin{itemize}
        \item $\{ 1 \}$ 
        \item $\{ r^{2k + 1} \}$
        \item $\{ r^t, r^{4k + 2 - t} \}$ for $1 \le t \le 2k$
        \item $\{ sr^{2t} \mid 0 \le t \le 2k \}$
        \item $\{ sr^{2t+1} \mid 0 \le t \le 2k \}.$
    \end{itemize}

    One conjugacy class of reflections corresponds to those in which the reflection line passes through vertices, and the other corresponds to those in which it passes through edges (i.e., the blue and red lines in \cref{fig:dihedral_geometric_representation}).
\end{recollection}

We now compute the protection group for $\D_{4k + 2}$. The proof starts with similar techniques, but becomes more sophisticated as we shall see.

\begin{theorem} 
    \label{thm:prot_for_d_4k+2}
    $\Prot(\D_{4k + 2}) \cong (\Zb/2)^{2k}$.
\end{theorem}

We prove the theorem by constructing a homomorphism $\varphi \colon \Prot(\D_{4k + 2}) \to(\Zb/2)^{2k}$, and proving this homomorphism is one-to-one and onto.

\begin{definition}
    For $1 \le t \le 2k$, we say a component is of \emph{type~$t$} if its arcs are colored by elements in $\{r^t, r^{n - t}\}$. Furthermore, we say a component is of type \emph{primary} if its arcs are colored by elements in $\{sr^{2t}\}$, and \emph{secondary} if its arcs are colored by elements in $\{sr^{2t + 1}\}$.
\end{definition}

\begin{definition}
    Define $\varphi \colon \Prot(\D_{4k + 2}) \to (\Zb/2)^{2k}$ such that for each $\D_{4k + 2}$-colored link $L$, the $i$-th component of $\varphi(L)$ is the number of times a strand of type~$i$ passes under a primary strand, modulo $2$.
\end{definition}

\begin{figure}[htbp]
    \centering
    \def\svgwidth{\textwidth}
\begingroup%
  \makeatletter%
  \providecommand\color[2][]{%
    \errmessage{(Inkscape) Color is used for the text in Inkscape, but the package 'color.sty' is not loaded}%
    \renewcommand\color[2][]{}%
  }%
  \providecommand\transparent[1]{%
    \errmessage{(Inkscape) Transparency is used (non-zero) for the text in Inkscape, but the package 'transparent.sty' is not loaded}%
    \renewcommand\transparent[1]{}%
  }%
  \providecommand\rotatebox[2]{#2}%
  \newcommand*\fsize{\dimexpr\f@size pt\relax}%
  \newcommand*\lineheight[1]{\fontsize{\fsize}{#1\fsize}\selectfont}%
  \ifx\svgwidth\undefined%
    \setlength{\unitlength}{283.46456693bp}%
    \ifx\svgscale\undefined%
      \relax%
    \else%
      \setlength{\unitlength}{\unitlength * \real{\svgscale}}%
    \fi%
  \else%
    \setlength{\unitlength}{\svgwidth}%
  \fi%
  \global\let\svgwidth\undefined%
  \global\let\svgscale\undefined%
  \makeatother%
  \begin{picture}(1,0.73)%
    \lineheight{1}%
    \setlength\tabcolsep{0pt}%
    \put(0,0){\includegraphics[width=\unitlength,page=1]{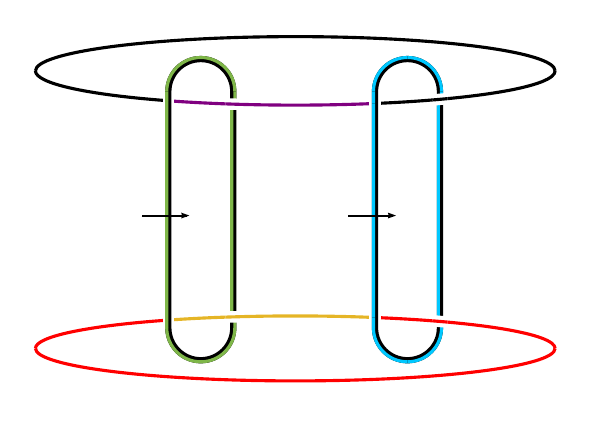}}%
    \put(0.49995695,0.56216717){\color[rgb]{0,0,0}\makebox(0,0)[t]{\lineheight{1.25}\smash{\begin{tabular}[t]{c}$sr^2$\end{tabular}}}}%
    \put(0.49995695,0.67891717){\color[rgb]{0,0,0}\makebox(0,0)[t]{\lineheight{1.25}\smash{\begin{tabular}[t]{c}$s$\end{tabular}}}}%
    \put(0.31943238,0.3314786){\color[rgb]{0,0,0}\makebox(0,0)[t]{\lineheight{1.25}\smash{\begin{tabular}[t]{c}$r$\end{tabular}}}}%
    \put(0.68662861,0.3314786){\color[rgb]{0,0,0}\makebox(0,0)[t]{\lineheight{1.25}\smash{\begin{tabular}[t]{c}$r^2$\end{tabular}}}}%
  \end{picture}%
\endgroup%

    \caption{$\varphi(L) = (1, 1)$, as both types pass under the primary strand once.}
    \label{fig:phi_calculation_example}
\end{figure}

Refer to \cref{fig:phi_calculation_example} for an example.


\subsubsection*{\texorpdfstring{$\varphi$}{phi} is onto}
\begin{figure}[htbp]
    \centering
    \def\svgwidth{0.7\textwidth}
\begingroup%
  \makeatletter%
  \providecommand\color[2][]{%
    \errmessage{(Inkscape) Color is used for the text in Inkscape, but the package 'color.sty' is not loaded}%
    \renewcommand\color[2][]{}%
  }%
  \providecommand\transparent[1]{%
    \errmessage{(Inkscape) Transparency is used (non-zero) for the text in Inkscape, but the package 'transparent.sty' is not loaded}%
    \renewcommand\transparent[1]{}%
  }%
  \providecommand\rotatebox[2]{#2}%
  \newcommand*\fsize{\dimexpr\f@size pt\relax}%
  \newcommand*\lineheight[1]{\fontsize{\fsize}{#1\fsize}\selectfont}%
  \ifx\svgwidth\undefined%
    \setlength{\unitlength}{283.46456693bp}%
    \ifx\svgscale\undefined%
      \relax%
    \else%
      \setlength{\unitlength}{\unitlength * \real{\svgscale}}%
    \fi%
  \else%
    \setlength{\unitlength}{\svgwidth}%
  \fi%
  \global\let\svgwidth\undefined%
  \global\let\svgscale\undefined%
  \makeatother%
  \begin{picture}(1,0.74)%
    \lineheight{1}%
    \setlength\tabcolsep{0pt}%
    \put(0,0){\includegraphics[width=\unitlength,page=1]{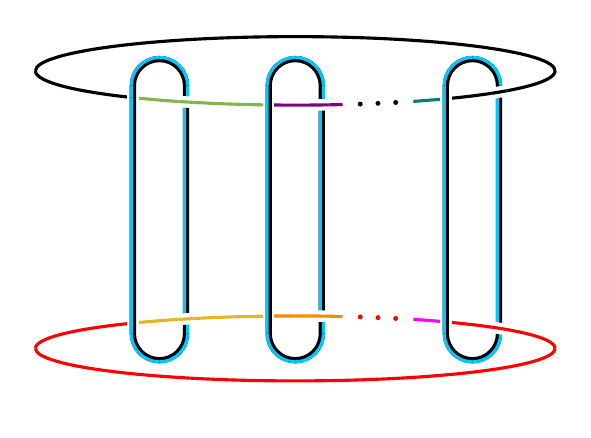}}%
    \put(0.49995693,0.69171355){\color[rgb]{0,0,0}\makebox(0,0)[t]{\lineheight{1.25}\smash{\begin{tabular}[t]{c}$s$\end{tabular}}}}%
    \put(0.49995693,0.0474886){\color[rgb]{0,0,0}\makebox(0,0)[t]{\lineheight{1.25}\smash{\begin{tabular}[t]{c}$sr^{2k+1}$\end{tabular}}}}%
    \put(0.19243403,0.3586829){\color[rgb]{0,0,0}\makebox(0,0)[t]{\lineheight{1.25}\smash{\begin{tabular}[t]{c}$r^t$\end{tabular}}}}%
  \end{picture}%
\endgroup%

    \caption{A link with the desired value under $\varphi$}
    \label{fig:generator_yielder_example}
\end{figure}

\begin{lemma}
    For each $1 \le i \le 2k$, there is a link $L_i$ such that $\varphi(L_i)$ 
    is zero in every component except the $i$-th.
\end{lemma}
\begin{proof}
    Consider the following example with one primary component, one secondary component, and $2k + 1$ components of type~$i$ (we denote this link by $L_i$), as in \cref{fig:generator_yielder_example}. Note that this is a valid coloring as $s$ and $sr^{2k + 1}$ return to their initial colors when conjugated $2k + 1$ times by $r^i$. The value of $\varphi(L_i)$ is as desired, as there are $2k + 1$ undercrossings with the primary component, and $2 \nmid 2k + 1$.

\end{proof}

\begin{theorem}
    $\varphi \colon \Prot(\D_{4k + 2}) \to (\Zb/2)^{2k}$ is onto.
\end{theorem}
\begin{proof}
    The elements $\varphi(L_1), \varphi(L_2), \dots, \varphi(L_{2k})$ generate $(\Zb/2)^{2k}$ and are in the image of $\varphi$. Therefore, $\varphi$ is onto.
\end{proof}

\subsubsection*{\texorpdfstring{$\varphi$}{phi} is one-to-one}
It suffices to prove that the kernel is trivial. Suppose $L$ is an element in the kernel. We proceed by simplifying the link, as we did for $\D_{2k + 1}$.

\begin{lemma}
    The link $L$ can evolve so that at most one primary and one secondary component exist.
\end{lemma}
\begin{proof}
    The proof is similar to that of \cref{lmm:joining_primary_components}.
\end{proof}

\begin{lemma}
    The link $L$ can evolve so that the primary and secondary components, if present, are unknots.
\end{lemma}
\begin{proof}
    The proof is similar to that of \cref{lmm:trivializing_the_primary}.
\end{proof}

If neither the primary nor secondary component exists, the problem is trivial as the remaining components commute. Furthermore, if only one of the components exists, one can proceed similarly and detach the components by surgery. Therefore, one can assume that both components exist.

One can trivialize the link further.

\begin{lemma}
    The link $L$ can evolve so that the secondary component only crosses over the primary component.
\end{lemma}
\begin{proof}
    The core idea is that for each primary color, there is one secondary color which commutes with it. Therefore, if a secondary strand of color $sr^{2t' + 1}$ passes under a strand of color $sr^{2t}$, it can be conjugated to $sr^{2t + (2k + 1)}$ by a technique similar to switching, and then crossed over $sr^{2t}$.
\end{proof}

With the help of this lemma, one can separate the two components entirely.

For simplicity, we remove the components colored by $r^0$ and $r^{2k + 1}$, which commute with everything.

We separate the primary and secondary components for an easier representation of the link. This will give a simpler and more intuitive understanding of how other components wind around these two.

\begin{figure}[htbp]
    \centering
    \includegraphics[width=\textwidth]{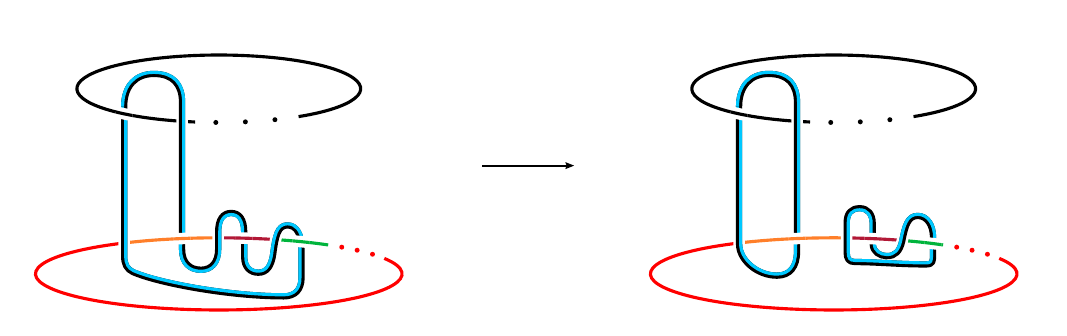}
    \caption{An example of a component broken down into two pieces.}
    \label{fig:breaking_down_the_component}
\end{figure}

Choose a component of type $t \in \{1, 2, \dots, 2k\}$. Equip it with a base point and an orientation. The component winds around the primary and secondary loops evenly many times in total, as the coloring must be consistent. If it winds around the two special components (primary or secondary) $2m$ times in total, it means that it can be broken down into $m$ pieces, each piece winding around the special components twice (\cref{fig:breaking_down_the_component}).

\begin{figure}[htbp]
    \centering

    \begin{subfigure}{0.9\textwidth}
        \centering
        \labeledfig{a} {%
            \includegraphics[width=\linewidth]{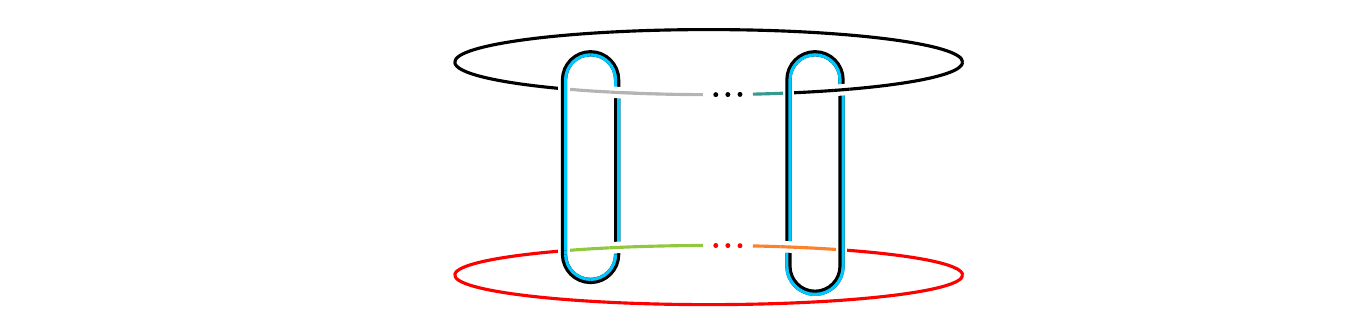}%
        }%
        \phantomcaption
        \label{fig:regular_vs_twisted}
    \end{subfigure}
    
    \hfill

    \begin{subfigure}{0.9\textwidth}
        \centering
        \labeledfig{b} {%
            \includegraphics[width=\linewidth]{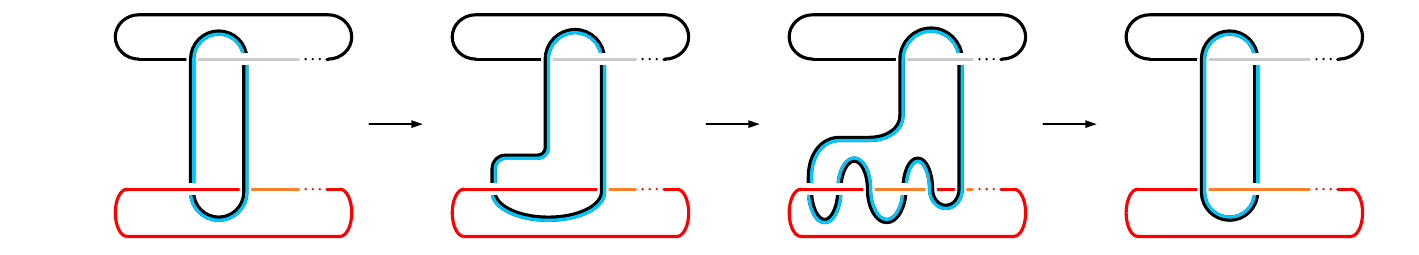}
        }
        \phantomcaption
        \label{fig:twisted_to_regular}
    \end{subfigure}

    \caption{(a) A component in the \emph{regular} form (left) and a component in the \emph{twisted} form (right). 
    (b) Transforming a twisted component into a regular component (the first step is performed via surgery).}
\end{figure}

If a piece winds around the same special loop twice, it is either an isolated unknot immediately or can be separated by surgery. The remaining components wind around both the primary and the secondary component. These components take one of the two forms presented in \cref{fig:regular_vs_twisted}. However, the \emph{twisted} form can be transformed into the \emph{regular} form as in \cref{fig:twisted_to_regular}. 

We have thus established the following simplification lemma:
\begin{lemma}
    $L$ can be transformed so that each component of type $1, 2, \dots, 2k$ 
    winds around each of the primary and secondary components once in the regular form.
\end{lemma}

It is also clear that the $i$-th component of $\varphi(L)$ is the number of type~$i$ components, modulo $2$.

We simplify the link by keeping the primary component at the top and the secondary component at the bottom as rectangles, with other components being chain-like rings connecting the lower side of the primary rectangle to the upper side of the secondary rectangle.

A component is identified by its topmost color and is called a \emph{standard loop of type~$t$}, where $t \in \{1, 2, \dots, 2k\} \cup \{2k + 2, 2k + 3, \dots, 4k + 1\}$.

\begin{figure}[htbp]
    \centering
    \def\svgwidth{0.5\textwidth}
\begingroup%
  \makeatletter%
  \providecommand\color[2][]{%
    \errmessage{(Inkscape) Color is used for the text in Inkscape, but the package 'color.sty' is not loaded}%
    \renewcommand\color[2][]{}%
  }%
  \providecommand\transparent[1]{%
    \errmessage{(Inkscape) Transparency is used (non-zero) for the text in Inkscape, but the package 'transparent.sty' is not loaded}%
    \renewcommand\transparent[1]{}%
  }%
  \providecommand\rotatebox[2]{#2}%
  \newcommand*\fsize{\dimexpr\f@size pt\relax}%
  \newcommand*\lineheight[1]{\fontsize{\fsize}{#1\fsize}\selectfont}%
  \ifx\svgwidth\undefined%
    \setlength{\unitlength}{198.42519685bp}%
    \ifx\svgscale\undefined%
      \relax%
    \else%
      \setlength{\unitlength}{\unitlength * \real{\svgscale}}%
    \fi%
  \else%
    \setlength{\unitlength}{\svgwidth}%
  \fi%
  \global\let\svgwidth\undefined%
  \global\let\svgscale\undefined%
  \makeatother%
  \begin{picture}(1,0.92857143)%
    \lineheight{1}%
    \setlength\tabcolsep{0pt}%
    \put(0,0){\includegraphics[width=\unitlength,page=1]{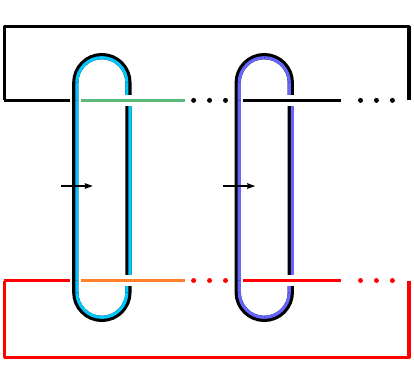}}%
    \put(0.23064625,0.42746739){\color[rgb]{0,0,0}\makebox(0,0)[t]{\lineheight{1.25}\smash{\begin{tabular}[t]{c}$r$\end{tabular}}}}%
    \put(0.64166944,0.42746739){\color[rgb]{0,0,0}\makebox(0,0)[t]{\lineheight{1.25}\smash{\begin{tabular}[t]{c}$r^2$\end{tabular}}}}%
  \end{picture}%
\endgroup%

    \caption{A type~$1$ standard loop on the left and a type~$2$ standard loop on the right.}
    \label{fig:chain_structure}
\end{figure}

This structure is explained in \cref{fig:chain_structure}.

One more simplification is required.

\begin{figure}[htbp]
    \centering
    \includegraphics[width=0.9\textwidth]{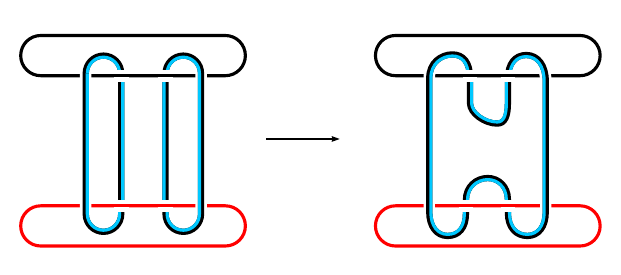}
    \caption{Canceling a standard loop of type~$t$ with one of color type~$4k + 2 - t$}
    \label{fig:canceling_standard_loops}
\end{figure}

\begin{lemma}
    The link $L$ can evolve so that standard loops of type~$t$ and type~$4k + 2 - t$ do not exist simultaneously.
    \label{lmm:canceling_standard_loops}
\end{lemma}
\begin{proof}
    The two standard loops can be brought adjacent to each other. By locally reconnecting them from the middle, they can be removed as an isolated unknot, as illustrated in \cref{fig:canceling_standard_loops}.
\end{proof}

\begin{figure}
    \centering
    \includegraphics[width=\textwidth]{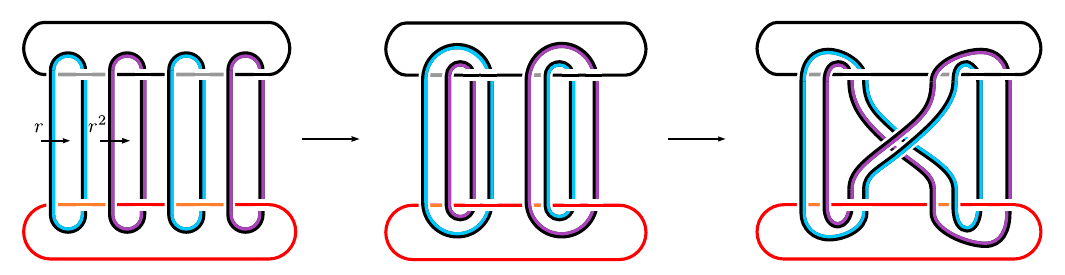}
    \caption{The two collections are attached to each other for $k = 1$.
    Note that the chains were slightly deformed to create this shape.}
    \label{fig:forming_the_simple_link}
\end{figure}

\begin{figure}
    \centering
    \includegraphics[width=\textwidth]{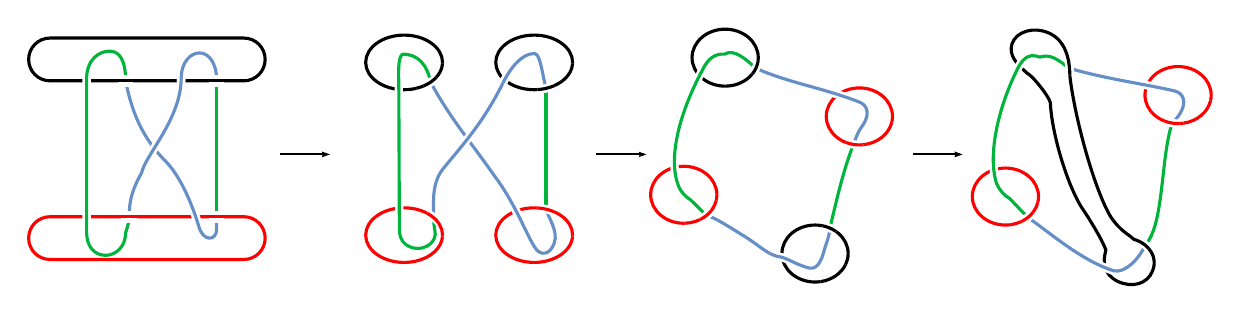}
    \caption{Deforming the strands to eliminate of the primary link.
    The green and blue strands represent the collection of numbered strands.}
    \label{fig:opening_the_simple_link}
\end{figure}

We now possess the machinery required to complete the proof of the theorem.
\begin{proof}[Proof of \cref{thm:prot_for_d_4k+2}]
    If $\varphi(L) = 0$, then we may assume using \cref{lmm:canceling_standard_loops}, for each $1 \le t \le 2k$,
    either evenly many standard loops exist of color $t$ and zero of color $4k + 2 - t$, or vice versa.
    Assume without loss of generality that evenly many exist of color $r^t$.
    Let the number of them be $2c_t$. 
    One can put $c_t$ on the left and $c_t$ on the right. 
    If this is done for all $t \in \{1, 2, \dots, 2k\}$, a similar collection of standard loops is concentrated on the left and on the right.
    Let the sum of the exponents for this collection be $\sigma$, and note that an arc of color $sr^\alpha$ is conjugated to $sr^{\alpha + 2\beta}$ if it passes under a collection of standard loops whose exponents sum to $\beta$. 
    Therefore, as $sr^{4\sigma} = s$ (the sum of the standard loops is $2\sigma$ in total), $sr^{2\sigma} = s$ ($4 \nmid 4k + 2$).
    Therefore, the special strands' colors return to their initial upon passing under the first collection.

    We sort the two collections in the opposite order, and in the two-dimensional diagram, arrange them such that they appear inside of one another. We then connect them as shown in \cref{fig:forming_the_simple_link}. 

    As \cref{fig:forming_the_simple_link} suggests, the collection can now be considered a single strand. 
    With some operations, one can deform the link (\cref{fig:opening_the_simple_link}) and remove the primary component.
    Similar steps can be taken for the secondary component, and the theorem is proven.
\end{proof}

\section{References}
\printbibliography[heading=none]
\end{document}